# A counterexample to Lagrangian Poincaré recurrence in dimension four

Joel Schmitz

**Abstract**

Counterexamples to Lagrangian Poincaré recurrence were recently found in dimensions greater than six by Broćić and Shelukhin. We construct counterexamples in dimension four using almost toric fibrations.

## 1 Introduction

Guided by the classical Poincaré recurrence for sets of positive measure and by the many Lagrangian intersection results in symplectic topology, Viterbo and Ginzburg around 2010 independently formulated the Lagrangian Poincaré recurrence conjecture (LPR):

> For every symplectic manifold $(M, \omega)$ it holds true that for every closed Lagrangian submanifold $L \subset M$ and every compactly supported Hamiltonian diffeomorphism $\varphi$ of $M$, there exists a sequence of integers $k_i \to \infty$ such that $\varphi^{k_i}(L) \cap L \neq \emptyset$.

We refer to [7] for a detailed discussion of this conjecture. In dimension $\geq 6$, counterexamples to LPR were recently constructed by Broćić and Shelukhin [2] for certain product tori in the ball, and hence in any symplectic manifold. They applied Chekanov's construction from [3]. In this note, we construct the first counterexamples in dimension four. Our ambient manifolds are 1-point blow-ups of $S^2 \times S^2$, and the Lagrangians are tori. The Hamiltonian diffeomorphisms without recurrence for these tori will be obtained with the help of Symington's structure theorem for almost toric symplectic 4-manifolds. The construction is similar to [9, Section 6], in that it is given by a loop of nodal slides. Our arguments also work in all closed symplectic toric 4-manifolds except for a short explicit list. See

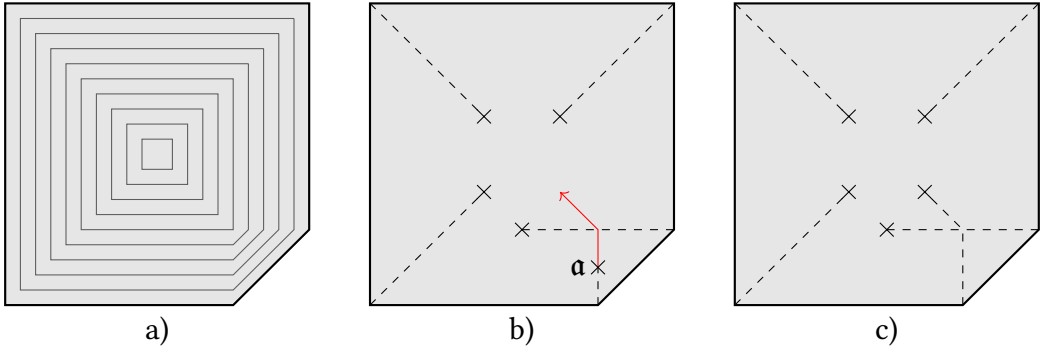

Figure 1: Constructing the almost toric fibration $\pi_0$.

Remarks 3.2 and 3.3 for a more detailed discussion. We also get counterexamples in dimension greater than 4, with quite different character than the counterexamples in [2]. See Remark 3.5

The counterexamples in [2] and here mean that LPR does not hold unconditionally. On the other hand, Lagrangian Poincaré recurrence was established in [8] for a certain class of Hamiltonian diffeomorphisms of $\mathbb{C}P^n$ called pseudo-rotations. This was generalized in [10]. Furthermore, the Lagrangian packing number of certain Lagrangian tori in non-monotone $S^2 \times S^2$ is finite, which in particular implies Lagrangian Poincaré recurrence, as shown by Polterovich–Shelukhin [15]. See also the previous work by Mak–Smith [13].

**Acknowledgement** We thank Felix Schlenk for fruitful inputs and discussions and Joé Brendel for his critical remarks.

## 2 The construction

Let $M$ be the one point blow-up of $(S^2 \times S^2, \omega_{a,b})$ by weight $c$, where $\omega_{a,b}$ is the product symplectic form giving area $a$ and $b$ to the factors with $a \geq b$ and $c < \frac{b}{2}$.

We will first construct a symplectomorphism $\varphi$ and a Lagrangian $L$ such that the Lagrangians $\{\varphi^n(L)\}_{n \in \mathbb{N}}$ are pairwise disjoint. Then we conclude using Lemma 2.2 that $\varphi$ (if $a \neq b$) or $\varphi^2$ (if $a = b$) must be a Hamiltonian symplectomorphism.

Let $\mu \colon M \to \Delta$ be the toric moment map where $\Delta$ is the Delzant polygon in Figure 1 a). $\Delta$ can be written as the non-negative locus of the function

$$\mathcal{F}_\Delta(x_1, x_2) = \min\left\{\frac{a}{2} - |x_1|, \frac{b}{2} - |x_2|, x_2 - x_1 + \frac{a+b}{2} - c\right\}.$$

$\mathcal{F}_\Delta$ measures the integral affine distance to the boundary of $\Delta$. The level sets of $\mathcal{F}_\Delta$ are also depicted in Figure 1 a). By making a nodal trade ([16, Theorem 6.5]) in each corner we obtain an almost toric fibration $\pi'_0 \colon X \to B'_0$, whose base diagram is depicted in Figure 1 b). The solid red broken line segment emanating from node $\mathfrak{a}$ in Figure 1 b) is a straight line segment in the integral affine structure on $B'_0$ inherited from the almost toric fibration and is contained in the eigenline of node $\mathfrak{a}$. Modifying $\pi'_0$ by a nodal slide ([16, Proposition 6.2]) along this line segment, we obtain the almost toric fibration $\pi_0 : M \to B_0$ depicted in Figure 1 c).

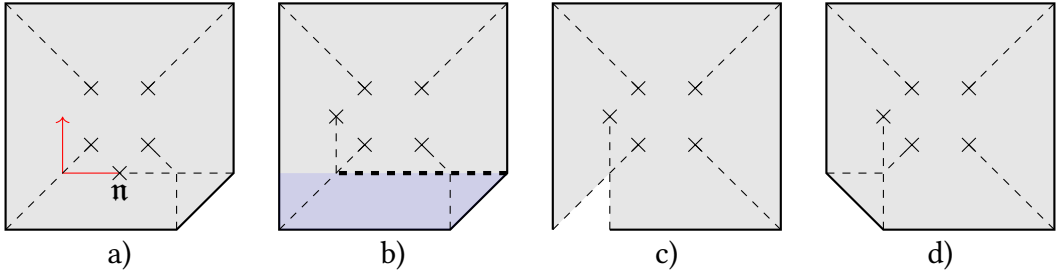

Figure 2: Nodal slide and change of branch cut

Note that $\max_{x \in \Delta}\{\mathscr{F}_\Delta(x)\} = \frac{b}{2}$. Define the $\mathbb{R}$-action $t\cdot$ on $\Delta \setminus \mathscr{F}_\Delta^{-1}(\frac{b}{2})$ such that $t \cdot x$ is the translation of $x$ by integral affine distance $t$ along the level set $\mathscr{F}_\Delta^{-1}(\mathscr{F}_\Delta(x))$ in the counter-clockwise direction.

Given $\varepsilon \in (0, \min\{c, \frac{b}{2} - c\})$, we will construct $\varphi \in \mathrm{Symp}(M)$ such that $\varphi$ acts on the fibres of $\pi_0$ as follows: Let $h = \mathscr{F}_\Delta(x)$. If $h > c + \varepsilon$, then $\varphi(\pi_0^{-1}(x)) = \pi_0^{-1}(x)$. If $h < c - \varepsilon$, then $\varphi(\pi_0^{-1}(x)) = \pi_0^{-1}((c - h) \cdot x)$. The integral affine length of a level set $\mathscr{F}_\Delta^{-1}(h)$ with $0 \leq h < c$ is $2(a+b) - c - 7h$. So for any $h \in (0, c - \varepsilon)$ with $\frac{c-h}{2(a+b)-c-7h}$ irrational and $x \in \mathscr{F}_\Delta^{-1}(h)$ the orbit of $\pi_0^{-1}(x)$ under $\varphi$ is dense in the level set $\mathscr{F}_\Delta^{-1}(h)$ and the Lagrangians $\{\varphi^n(\pi_0^{-1}(x))\}_{n \in \mathbb{N}}$ are pairwise disjoint.

**Constructing $\varphi$.** Note that the node $\mathfrak{n}$ in Figure 2 a) has $\mathscr{F}_\Delta(\mathfrak{n}) = c$. We modify the fibration $\pi_0$ by a sequence of four nodal slides, illustrated in Figure 2: First slide the node $\mathfrak{n}$ along the marked path in Figure 2 a) to obtain an almost toric fibration with base diagram Figure 2 b). During the nodal slide, using [5, Proof of Theorem 8.10] we may assume that we do not modify the fibration outside of $\mathscr{F}_\Delta^{-1}(c - \varepsilon, c + \varepsilon)$. Note that the branch cuts in the base diagram Figure 2 b) cut the base $B_0$ into three connected components. Removing the horizontal branch cut drawn thickly in Figure 2 b), we get the base diagram Figure 2 c), which has only one connected component. This diagram has a discontinuity in the lower left corner, which can be removed by introducing a branch cut and applying a suitable shear map (given below), giving Figure 2 d). The effect of the two steps b)→c)→d) is to apply the shear map

$$(x_1, x_2) \mapsto \left(x_1 + \left(c - \frac{b}{2} - x_2\right), x_2\right) \tag{1}$$

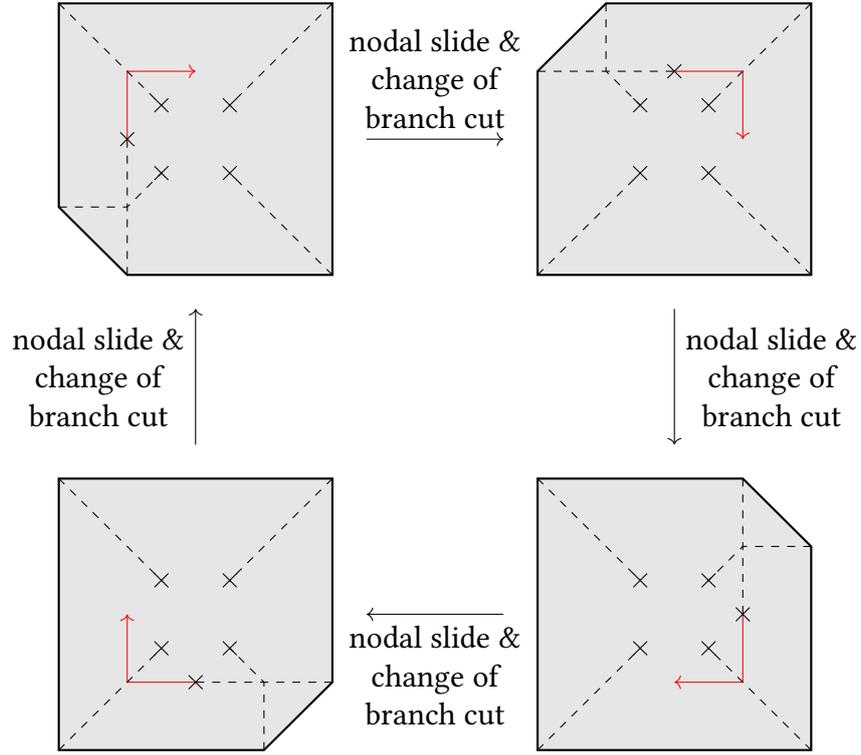

Figure 3: Constructing the almost toric fibration $\pi_1$.

to points with $x_2 < c - \frac{b}{2}$ (the region shaded blue in Figure 2 b) ), and the identity everywhere else, giving the base diagram Figure 2 d). Note that the diagram d) is identical to the diagram a) by a 90°-rotation. Repeat this process 3 more times as illustrated in Figure 3 until $\mathfrak{n}$ has returned to its original position. By this process we obtain an almost toric fibration $\pi_1 : M \to B_1$ such that, for $x \in \Delta$, if $h > c + \varepsilon$ the fibres $\pi_1^{-1}(x)$ and $\pi_0^{-1}(x)$ are equal and if $h < c - \varepsilon$, by combining the effect of the four shear maps (1), we have

$$\pi_1^{-1}(x) = \pi_0^{-1}((c - h) \cdot x) . \tag{2}$$

Now [16, Corollary 5.4] and [5, Proof of Theorem 8.5] gives a symplectomorphism $\varphi \in \text{Symp}(M)$ such that outside a small neighbourhood of the nodes

$$\varphi(\pi_0^{-1}(x)) = \pi_1^{-1}(x) .$$

From (2) we conclude that for $h = \mathscr{F}_\Delta(x) < c - \varepsilon$

$$\varphi(\pi_0^{-1}(x)) = \pi_1^{-1}(x) = \pi_0^{-1}((c - h) \cdot x) .$$

*Remark* 2.1. As pointed out by the referee, here is an equivalent manner of seeing the effect of the nodal slide in Figure 3: Figure 4 a) gives an alternative base diagram of the almost toric base $B_0$ pictured in Figure 1 c), which can also be obtained from $(S^2 \times S^2, \omega_{a,b})$ by a non-toric blow-up as in [5, Section 9.1]. Similarly to [5, Section 10.1], we may replace the four branch cuts emanating from the nodes $\mathfrak{a}, \mathfrak{b}, \mathfrak{c}, \mathfrak{d}$ by the tetrapod branch cut $T$ highlighted in Figure 4 a). Denote the branch cut emanating from the node $\mathfrak{n}$ by $\ell$. Then $B_0 \setminus (T \cup \ell)$ is homotopy equivalent to $S^1$, and a base diagram of its universal cover is given by Figure 4 b). The concatenation of all the nodal slides in Figure 3 can be written as one nodal slide marked in Figure 4 a) and b) by the solid red line emanating from $\mathfrak{n}$. Performing the nodal slide, we get the base diagram Figure 4 c) and the almost toric fibration $\pi_1 : M \to B_1$ as above. As described in [5, Section 7.2] we can rotate the branch cut $\ell$ emanating from $\mathfrak{n}$ until it matches the branch cut in diagram Figure 4 b) to obtain the base diagram Figure 4 d) for the universal cover of $B_1 \setminus (T \cup \ell)$. The combined effect of the steps b)→c)→d) is to apply the shear matrix $A = \begin{pmatrix} 1 & 1 \\ 0 & 1 \end{pmatrix}$ to points below the eigenline of $\mathfrak{n}$: The "missing" triangle below $\mathfrak{n}$ is spanned by vectors in directions $(a, -1), (a+1, -1)$ for some value of $a \in \mathbb{R}$, so when a point $x$ touches the sliding or rotating branch cut, in the base diagram the point jumps over the gap given by the missing triangle according the shear matrix $A$. As above we use [5, Proof of Theorem 8.10] to ensure the nodal slide only modifies fibres near the eigenline of $\mathfrak{n}$, and [16, Corollary 5.4] and [5, Proof of Theorem 8.5] to get the desired symplectomorphism.

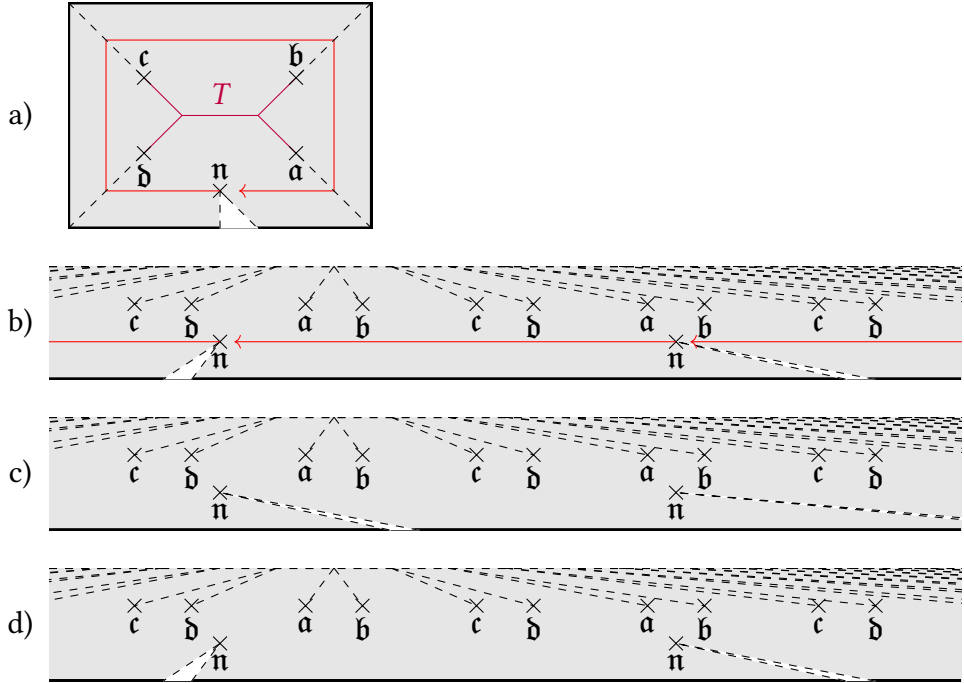

Figure 4: An alternative picture of the nodal slides used in the construction of the symplectomorphism $\varphi$.

We conclude the construction by showing that $\varphi$ or $\varphi^2$ is a Hamiltonian diffeomorphism.

**Lemma 2.2.** *Let $\psi$ be a symplectomorphism of M.*
  1. *If $a \neq b$, then $\psi$ is a Hamiltonian diffeomorphism.*
  2. *If $a = b$, then $\psi^2$ is a Hamiltonian diffeomorphism.*

*Proof.* It is shown in [12, Cor. 1.3] that the symplectic mapping class group of $M$ is a finite reflection group generated by Dehn-twists along embedded Lagrangian spheres, and that a generating set is given by the classes $\xi \in H_2(M;\mathbb{Z})$ which can be represented by smoothly embedded spheres and satisfy

$$\xi \cdot \xi = -2, \quad c_1(M)(\xi) = 0, \quad [\omega](\xi) = 0$$

Already the conditions $\xi \cdot \xi = -2$ and $c_1(M)(\xi) = 0$ imply that $\xi$ is the anti-diagonal class $-D$, see e.g. [6, Lemma 4.1.4]. Since the symplectic area of $-D$ vanishes only for $a = b$, and since in this case the Dehn-twist along the anti-diagonal has order at most 2, the symplectic mapping class group is trivial in case $a \neq b$ and has at most two elements in case $a = b$. Finally note that the symplectomorphisms in the component of the identity are Hamiltonian, since the first homology of $M$ vanishes. □

Alternatively [11, Corollary 4.6] yields a shorter but slightly weaker argument. See Remark 3.3.

## 3 Discussion

*Remark 3.1.* The set of $h \in (0, c - \varepsilon)$ such that for $x \in \mathscr{F}_\Delta^{-1}(h)$ the Lagrangians $\{\varphi^n(\pi_0^{-1}(x))\}_{n\in\mathbb{N}}$ are pairwise disjoint is dense and has full measure in $(0, c - \varepsilon)$.

*Remark 3.2.* The construction in Section 2 of a symplectomorphism $\varphi$ can be generalized to more toric symplectic 4-manifolds: Let $\Delta$ be a compact Delzant polygon and $\mu\colon X \to \Delta$ the toric symplectic manifold corresponding to $\Delta$. Denote by $\mathscr{F}_\Delta\colon \Delta \to \mathbb{R}$ the integral affine distance to the boundary of $\Delta$. Then we can repeat the construction in Section 2 if $\Delta$ has an edge of length $c < \max_{x\in\Delta}\{\mathscr{F}_\Delta(x)\}$ corresponding to a symplectic sphere of self-intersection $-1$. By [14], these are all closed symplectic toric 4-manifolds except:
  1. The five monotone symplectic toric 4-manifolds.
  2. $(S^2 \times S^2, \omega_{a,b})$ with $a \neq b$.
  3. The one or two-fold blow-up of $(S^2 \times S^2, \omega_{a,b})$ of size $c = \frac{\min\{a,b\}}{2}$.

*Remark 3.3.* Let $(M, \omega)$ be a closed toric symplectic 4-manifold. By [11, Corollary 4.6], the symplectic mapping class group of $(M, \omega)$ is finite.

Thus the generalized construction in Remark 3.2 yields a symplectomorphism $\varphi$ such that for many fibre tori $L$, the tori $\{\varphi^n(L)\}_{n\in\mathbb{N}}$ are pairwise disjoint. Since the mapping class group is finite, $\psi = \varphi^k$ is Hamiltonian for some $k \geq 1$ and the tori $\{\psi^n(L)\}_{n\in\mathbb{N}}$ are pairwise disjoint. In particular, LPR does not hold for these toric symplectic 4-manifolds.

*Remark* 3.4. While our construction does not yield a result for $S^2 \times S^2$, in [1] Brendel and Kim show that most toric fibres of non-monotone $S^2 \times S^2$ have infinite Lagrangian packing number.

*Remark* 3.5. Taking the product of our counterexamples $(M, L, \varphi)$ to LPR with any pair $(M', L')$ and the identity of $M'$, we obtain many counterexamples to LPR in dimensions $\geq 6$. By our construction, every Lagrangian $\varphi^k(L) \times L'$ is an accumulation point of the orbit $\varphi^n(L) \times L'$. These counterexamples are therefore very different from the ones in [2], where the sequence of tori $\phi^n(L)$ converges to one Lagrangian torus.

Using Delzants construction [4], we can obtain any closed symplectic toric manifold by symplectic reduction from $\mathbb{C}^N$ where $N$ is the number of edges of the moment polytope. Lifting the Hamiltonians constructed above, we also get counterexamples to LPR in symplectic balls in dimension $\geq 8$ (e.g. a non-monotone one point blow-up of $\mathbb{C}P^2$ can be obtained from $\mathbb{C}^4$ by symplectic reduction), and thus any symplectic manifold of dimension $\geq 8$.